\documentclass[11pt]{article}
\usepackage{amssymb,psfig}
\usepackage{amsmath}
\usepackage{times}\hfuzz=10pt 
\newtheorem{theorem}{Theorem}%[section]
\newtheorem{corollary}{Corollary}%[section]
\def\e{\varepsilon}

\def\defi{\stackrel{{\scriptscriptstyle \Delta}}{=}}

\def\a{\alpha}
\def\d{\delta}
\def\o{\omega}
\def\O{\Omega}

\def\F{{\cal F}}

\def\Arg{{\rm Arg\,}}

\def\Re{{\rm Re\,}}

\def\R{{\bf R}}

\def\b{\beta}

\def\C{{\bf C}}

\def\ww{\widetilde}

\def\oo{\bar}

\newcommand{\be}{\begin{equation}}
\newcommand{\ee}{\end{equation}}
\newcommand{\bd}{\begin{displaymath}}
\newcommand{\ed}{\end{displaymath}}
\newcommand{\ba}{\begin{array}{ll}}
\newcommand{\ea}{\end{array}}
\newcommand{\baa}{\begin{eqnarray}}
\newcommand{\eaa}{\end{eqnarray}}
\newcommand{\baaa}{\begin{eqnarray*}}
\newcommand{\eaaa}{\end{eqnarray*}}
\font\sm=cmr10
\def\oo{\bar}

\def\a{\alpha}

\def\hh{h}
\def\HH{H}
\def\a{a}\def\b{b}
\date{ Submitted: December 29, 2010. Revised:
June 23, 2011}
\title{
On sub-ideal causal smoothing filters
 \footnote{Accepted to {\em Signal Processing}: Submitted: December 29, 2010. Revised:
June 23, 2011}}
\author{
Nikolai Dokuchaev\\  {\sm Department of Mathematics \& Statistics,
Curtin University,}\\
{\sm  GPO Box U1987, Perth, 6845 Western Australia}\\ {\sm email
N.Dokuchaev@curtin.edu.au. Tel.: 61 8 92663144.}}
 
\begin{document}
 \vspace{-0.5cm}      \maketitle
\begin{abstract}  Smoothing causal linear time-invariant filters are studied  for
continuous time processes. The paper suggests a family of  causal
filters with almost exponential damping of the energy on the higher
frequencies.  These filters are sub-ideal meaning that a faster
decay of the frequency response would lead to the loss of causality.
%\\
%PACS 2008 numbers: 02.30.Mv, %    Approximations and expansions
%02.30.Nw,  %  Fourier analysis
%02.30.Yy, %    Control theory
%07.05.Mh,  %  Neural networks, fuzzy logic, artificial intelligence
%07.05.Kf %Data analysis: algorithms and implementation; data management
\end{abstract}
\section{Introduction}  The paper studies   smoothing filters  for
continuous time processes.  The consideration is restricted by the
causal continuous time linear time-invariant filters (LTI filters),
i.e. linear  filters represented as convolution integrals over the
historical data. These filters are used in dynamic smoothing, when
the future values of the process are not available.
\par
In the frequency domain, smoothing means reduction of the energy on
the higher frequencies. In particular, an ideal low-pass filter is a
smoothing filter. However, this filter is not causal, i.e., it
requires the future value of the process. Moreover,  a filter with
exponential decay of the frequency response also cannot be causal
\cite{D10}. It follows from the fact  that a sufficient rate of
decay of energy  on higher frequencies implies some predictability
of the processes; on the other hand, a causal filter cannot
transform a general kind of a process into a predictable process.
The classical result is Nyquist-Shannon-Kotelnikov interpolation
theorem that implies that if a process is band-limited then it is
predictable (see, e.g., \cite{rema}-\cite{Br},\cite{D08},
\cite{Higgins}-\cite{M}, \cite{P}-\cite{W}). Recently, it was found
 that processes with exponential decay of energy on the
higher frequencies are weakly predictable on a finite time horizon
\cite{D10}.
\par
We suggest a family of causal smoothing filters  with "almost"
exponential rate of damping the energy on the higher frequencies and
with the frequency response that can be selected to approximate the
real unity uniformly on an arbitrarily large interval. These filters
are sub-ideal in the sense that their effectiveness in the damping
of higher frequencies cannot be exceeded; a faster decay of the
frequency response is not possible for causal filters. This is
because  this family of causal filters approximates the exponential
decay rate of a reference set of non-causal filters (\ref{NC}).
\section{Problem setting}
Let $x(t)$ be a continuous-time process, $t\in\R$. The output of a
linear filter is the process \baaa y(t)=\int_{-\infty}^{\infty}
\hh(t-\tau)x(\tau)d\tau,\eaaa where $\hh:\R\to\R$ is a given impulse
response function.

If $\hh(t)=0$ for $t<0$, then the output of the corresponding filter
is \baaa y(t)=\int_{-\infty}^{t} \hh(t-\tau)x(\tau)d\tau. \eaaa In
this case, the filter and the impulse response function are said to
be causal. The output of a causal  filter at time $t$ can be
calculated using only past historical values $x(\tau)|_{\tau\le t}$
of the currently observable continuous-time input process.

The goal is to approximate $x(t)$ by a smooth filtered process
$y(t)$ via selection of an appropriate causal impulse response
function $\hh(\cdot)$.
\par
We are looking for families of the causal smoothing impulse response
functions $\hh(\cdot)$  satisfying the following conditions.
\begin{itemize}
\item[(A)] The  outputs $y(\cdot)$ approximate processes $x(\cdot)$; the arbitrarily close approximation
can be achieved by selection of an appropriate impulse response from
the family.
\item[(B)] For processes $x(\cdot)\in L_2(\R)$, the outputs $y(\cdot)$ are infinitely differentiable
functions. On the higher frequencies, the frequency response of the
filter is as small as possible, to achieve the most effective
damping of the energy on the higher frequencies of $x$.
\item[(C)] The effectiveness of this family in the damping of the
higher frequencies cannot be exceeded; any faster decay of the
frequency response  would lead to the loss of causality.
\item[(D)] The effectiveness of this family in the damping of the
higher frequencies  approximates the effectiveness of some reference
family of non-causal smoothing filters with a reasonably fast decay
of the frequency response.
\end{itemize}
\par
Note that it is not a trivial task to satisfy Conditions (C)-(D).
For instance, consider a family of low-pass filters with increasing
pass interval $[-\Delta,\Delta]$, where $\Delta >1$. Clearly, the
corresponding smoothed processes approximate the original process as
$\Delta\to +\infty$, i.e., Condition (A) is satisfied. However, the
distance of the set of these ideal low-pass filters from the set of
all causal filters is positive \cite{rema}.
\par
For $x(\cdot)\in  L_2(\R)$, we denote by $X=\F x$ the function
defined on $i\R$ as the Fourier transform of $x(\cdot)$;
$$X(i\o)=(\F x)(i\o)= \int_{-\infty}^{\infty}e^{-i\o t}x(t)dt,\quad
\o\in\R.$$ Here $i=\sqrt{-1}$. For $x(\cdot)\in L_2(\R)$, the
Fourier transform $X$ is defined as an element of $L_2(\R)$ (more
precisely, $X(i\cdot)\in L_2(\R)$).
\par
Consider a reference family of "ideal" smoothing filters with the
frequency response \baa M_\mu(i\o)=e^{-\mu |\o|}, \quad \mu>0.
\label{NC}\eaa For these filters, Condition (A) is satisfied  as
$\mu\to 0$, and Conditions (B) is satisfied for all $\mu>0$.
However, these filters are non-causal: for any $x(\cdot)\in
L_2(\R)$, the output processes of these filters are weakly
predictable at time $t$ on a finite horizon $[t,t+\mu)$ \cite{D10}.
\par
To satisfy Conditions (A)--(D), we consider a family of causal
filters with impulse responses
$\{\hh_\nu(\cdot)\}_{\nu=1}^{\infty}\subset L_2(\R)$ and with the
corresponding Fourier transforms $\HH_\nu(i\o)$, such  that the
following more special Conditions (a)-(d) are satisfied.
\begin{itemize}
\item[(a)]  {\em  Approximation of identity operator:}
\subitem(a1) For any $\O>0$, $\HH_\nu(i\o)\to 1$ as $\nu\to +\infty$
uniformly in $\o\in[-\O,\O]$. \subitem(a2) For any  $x(\cdot)\in
L_2(\R)$, \baaa \|y_{\nu}(\cdot)-x(\cdot)\|_{L_2(\R)}\to 0\quad
\hbox{as}\quad \nu\to \infty,\eaaa where $y_\nu$ is the output
process \baaa y_\nu(t)=\int_{-\infty}^{t}
\hh_\nu(t-\tau)x(\tau)d\tau.\eaaa
\item[(b)] {\em Smoothing property:} For every
$\nu>0$, there exists $\rho>0$ such that for any $n\ge 1$, \baaa
\int_{-\infty}^{\infty}e^{|\o|^\rho}|\HH_{\nu}(i\o)|^nd\o<+\infty.\eaaa
\item[(c)]  {\em Sub-ideal smoothing:} For any  $\d>1$, there exists $\nu>0$ such that
for any $\O>0$  \baa \int_{\{ \o:\
|\o|\ge\O\}}\frac{|\log|\HH_\nu(i\o)||^\d}{1+\o^2}d\o=+\infty.
\label{delta} \eaa
\item[(d)]  {\em Approximation  of non-causal filters
(\ref{NC}) with respect to the effectiveness in damping:} For any
$\e>0$ and $\mu>0$, there exists $\nu=\nu(\mu)>0$ such that \baaa
\bigl\||\HH_{\nu}(i\o)|-|M_\mu(i\o)|\bigr\|_{L_2(\R)}\le \e. \eaaa
\end{itemize}
\par
Let us show that Conditions (a)-(d) ensure that Conditions (A)-(D)
are satisfied, in a certain sense. Clearly, Condition (a) ensures
that  Condition (A) is satisfied.
\par
Further, by Condition (b),  for any $k>0$ and $\nu>0$, \baaa
\int_{-\infty}^{\infty}(1+|\o|^k)^4|\HH_{\nu}(i\o)|^4d\o<+\infty.\eaaa
Let $x(\cdot)\in L_2(\R)$, $X=\F X$, and
$Y_\nu(i\o)=\HH_\nu(i\o)X(i\o)$. By H\"older inequality, it follows
that \baaa \int_{-\infty}^{\infty}(1+|\o|^k)^2|Y_{\nu}(i\o)|^2d\o\le
\left(\int_{-\infty}^{\infty}(1+|\o|^k)^4|\HH_{\nu}(i\o)|^4d\o\right)^
{1/2}\|X\|_{L_2(\R)}^{1/2} <+\infty.\eaaa Hence $y_{\nu}(t)$ has
derivatives  in $L_2(\R)$ of any order, and, therefore, is
infinitely differentiable in the classical sense. Therefore,
Condition (b) ensures that Condition (B) is satisfied.
\par
Let us show that Condition (c) ensures that  Condition (C) is
satisfied. Let $\d>1$ be fixed, and  let $\nu=\nu(\d)$ be such that
(\ref{delta}) holds. Let us show that the filter with the frequency
response  $\ww\hh=\F^{-1}\ww\HH$ cannot be causal for  some "better"
frequency response  $\ww \HH(i\o)$ such that \baa |\ww
\HH(i\o)|=o(|\HH_\nu(i\o)|)\quad \hbox{as}\quad |\o|\to +\infty.
\eaa More precisely, we will show that $\ww \hh$ cannot be causal
with a stronger condition that there exists $\O>0$ such that \baa
|\log|\ww \HH(i\o)||\ge |\log|\HH_\nu(i\o)||^{\d},\quad |\o|\ge
\O.\label{d1}\eaa In particular, this condition implies that
$\log|\HH_\nu(i\o)|/\log|\ww\HH(i\o)|\to 0$ as $|\o|\to +\infty$.
\par
The desired fact that $\ww\hh$ cannot be causal  can be seen from
the following. By Paley and Wiener Theorem \cite{PW}, the Fourier
transform $\HH(i\o)$ of a causal impulse response $\hh(\cdot)\in
L_2(\R)$ has to be such that\baaa
\int_{-\infty}^{\infty}\frac{|\log|\HH(i\o)||}{1+\o^2}d\o<+\infty
\label{est}\eaaa (see, e.g., \cite{Paarmann}, p.35).  Since
$\nu=\nu(\d)$ is such that (\ref{delta}) holds, it follows from
(\ref{d1}) that $\ww\hh$ cannot be causal. Therefore, Condition (c)
ensures that Condition (C) is satisfied.
\par
Finally,  Condition (d) ensures that Condition (D) is satisfied,
since the effectiveness of smoothing is defined by the rate of
damping  of the higher frequencies.
\section{A family of sub-ideal smoothing filters}
  Let $\C^+\defi\{z\in\C:\ \Re z> 0\}$.
  Let us consider a set of transfer functions
\baa \HH_{\a,\b,q}(s)\defi e^{-\a (s+\b)^q},\quad s\in\C^+.
\label{K}\eaa Here $\a>0$, $\b>0$, and $q\in[\oo q,1)$, are rational
numbers, $\oo q\in(0,1)$ is a given number.
 We mean the branch of $(s+\b)^q$
such that its argument is $q\Arg(s+\b)$, where $\Arg z\in
(-\pi,\pi]$ denotes the principal value of the argument of $z\in\C$.
This set was introduced in \cite{D07} as an auxiliary tool for
solution of a parabolic equation in the frequency domain.
\par
Let us consider the set of all transfer functions (\ref{K}) with
rational numbers $\a>0$, $\b>0$, and $q\in[1/2,1)$. We assume that
this countable set is counted as a sequence
$\{\HH_\nu\}_{\nu=1}^{\infty}$ such that $\a\to 0$, $b\to 0$, $q\to
1$ as $\nu\to +\infty$.
\begin{theorem}\label{ThM} Conditions (a)-(d) are
satisfied  for  the family of filters defined by the transfer
functions $\{\HH_\nu\}_{\nu=1}^{\infty}$. (Therefore, Conditions
(A)-(D) are satisfied for this family).
\end{theorem}
\par
{\em Proof of Theorem \ref{ThM}.}
 Let  $H^r$ be the Hardy space of holomorphic on $\C^+$
functions $h(p)$ with finite norm
$\|h\|_{H^r}=\sup_{\rho>0}\|h(\rho+i\o)\|_{L_r(\R)}$,
$r\in[1,+\infty]$ (see, e.g., \cite{Du}).
\par
Clearly, the functions $\HH_\nu(p)$ are holomorphic in $\C^+$, and
\baa \ln|\HH_\nu(s)|=-\Re(\a (s+\b)^q)=-\a|s+\b|^q\cos [q\Arg
(s+\b)].\label{logK}\eaa  In addition, there exists $M=M(\b,q)>0$
such that $\cos[q\Arg(p+\b)]>M$ for all $s\in \C^+$. It follows that
\baa |\HH_\nu(s)|\le e^{-\a M|s+\b|^q}<1,\quad s\in\C^+.
\label{Kest}\eaa Hence $\HH_\nu\in H^r$ for all $r\in [1,+\infty]$.
By Paley-Wiener Theorem, the inverse Fourier transforms
$\hh_\nu=\F^{-1}\HH_\nu(i\o)$ are causal impulse responses, i.e.,
$\hh_\nu(t)=0$ for $t<0$ (see, e.g., \cite{Yosida}, p.163).

Let $x\in L_2(\R)$, $X=\F x$, and $Y_\nu=\HH_\nu X$.

Let us show that Condition (a) holds. Since $\a\to 0$ as
$\nu\to+\infty$, it follows that $ \HH_\nu(i\o)\to 1$ as $\nu\to
+\infty$ for any $\o$ and that Condition (a1) holds. By Condition
(a1), $ Y_\nu(i\o)\to X(i\o)$ as $\nu\to +\infty$ for all $\o\in\R$.
In addition, $|\HH_\nu(i\o)|\le 1$. Hence $|Y_\nu(i\o)-X(i\o)|\le
2|X(i\o)|$. We have that
$\|X(i\o)\|_{L_2(\R)}=\|x\|_{L_2(\R)}<+\infty$.
 By Lebesgue
Dominance Theorem, it follows that \baaa \left\|Y_\nu(i\o)-X(i\o)
\right\|_{L_2(\R)}\to 0\quad\hbox{as}\quad \nu\to +\infty.
 \eaaa
Therefore, Condition (a) holds.

Let us show that Condition (b) holds. By (\ref{Kest}), it follows
that
 \baa |\HH_\nu(i\o)|\le e^{-\a M|\o|^q},\quad
\o\in\R. \label{Kest2}\eaa
 Therefore, Condition (b) holds  with any $\rho< q$.

To see that Condition (c) holds, it suffices to observe that
(\ref{delta}) holds if $\d q\ge 1$, i.e., $q\ge 1/\d$.

Let us show  that Condition (d) holds.  We assume that the family of
transfer functions  $\HH_{\a,\b,q}(\cdot)$ is counted as a sequence
$\{\HH_\lambda\}_{\lambda=1}^{\infty}$ such that  $\b\to 0$, $q\to
1$ as $\lambda\to +\infty$, with $\a=\mu/\cos(q\pi/2)$. We have that
 $\cos[q\Arg(i\o+\b)]\ge \cos(q\pi/2)$ for all $b>0$, $\o\in\R$ and
 $\cos[q\Arg(i\o+\b)]\to \cos(q\pi/2)$ as $b\to 0$ for all $q\le 1$, $\o\in\R$.

By (\ref{logK}), $ |\HH_\lambda(i\o)|\to |M_\mu(i\o)|$ as
$\lambda\to +\infty$ for all $\o\in\R$.
 Further, we have that
 \baa -\ln|\HH_\lambda(i\o)|=a\Re
[(i\o+\b)^q]=\a|i\o+\b|^q\cos [q\Arg (i\o+\b)]\nonumber\\\ge
\a|i\o+\b|^q\cos (q\pi/2) \ge \a|\o|^q\cos
(q\pi/2).\label{logK1}\eaa Hence $ -\ln|\HH_\lambda(i\o)|\ge
\mu|\o|^q$ and $ |\HH_\lambda(i\o)|\le e^{- \mu|\o|^q}.$
 Since $q\ge \oo q>0$, we have from
(\ref{Kest2}) that, for some constants $c_1>0$ and $c_2>0$,
$|\HH_\lambda(i\o)|+|M_\a(i\o)| \le c_1 \exp\left(-c_2(|\o|^{\oo
q}+|\o|)\right)$ for all $\o$.
 By Lebesgue
Dominance Theorem, it follows that \baa
\bigl\||\HH_\lambda(i\o)|-|M_\mu(i\o)| \bigr\|_{L_2(\R)}\to
0\quad\hbox{as}\quad \lambda\to +\infty.
 \label{L2}\eaa  Hence Condition (d) holds. This completes
the proof of Theorem \ref{ThM}. $\Box$
\par
Note that the sequence $\HH_\lambda(s)$ introduced above does not
ensure approximation described by Condition (a1), since $\a\to
+\infty$ and $\HH_\lambda(0)\to 0$ as $\lambda\to +\infty$. On the
other hand, the sequence $\HH_\nu(s)$ does not ensure approximation
(\ref{L2}). The following corollary shows a way toward the
combination of these approximation properties.
\begin{corollary}\label{corr1} Let $q\in[\oo q,1)$ and
$b>0$  be given. Let a sequence
$\{\HH_{n}(\cdot)\}=\{\HH_{\a,\b,q}(\cdot)\}$ be selected such that
$a\to 0$. Then Condition (a1) holds for this sequence, and
$|\HH_n(i\o)|\le e^{- c|\o|^q}$ for all $\o$, where $c= \a\cos
(q\pi/2)$.
\end{corollary}
\par Proof of Corollary \ref{corr1} follows from (\ref{logK1}).
\section{Illustrative examples} 
\begin{figure}[ht]
\centerline{\psfig{figure=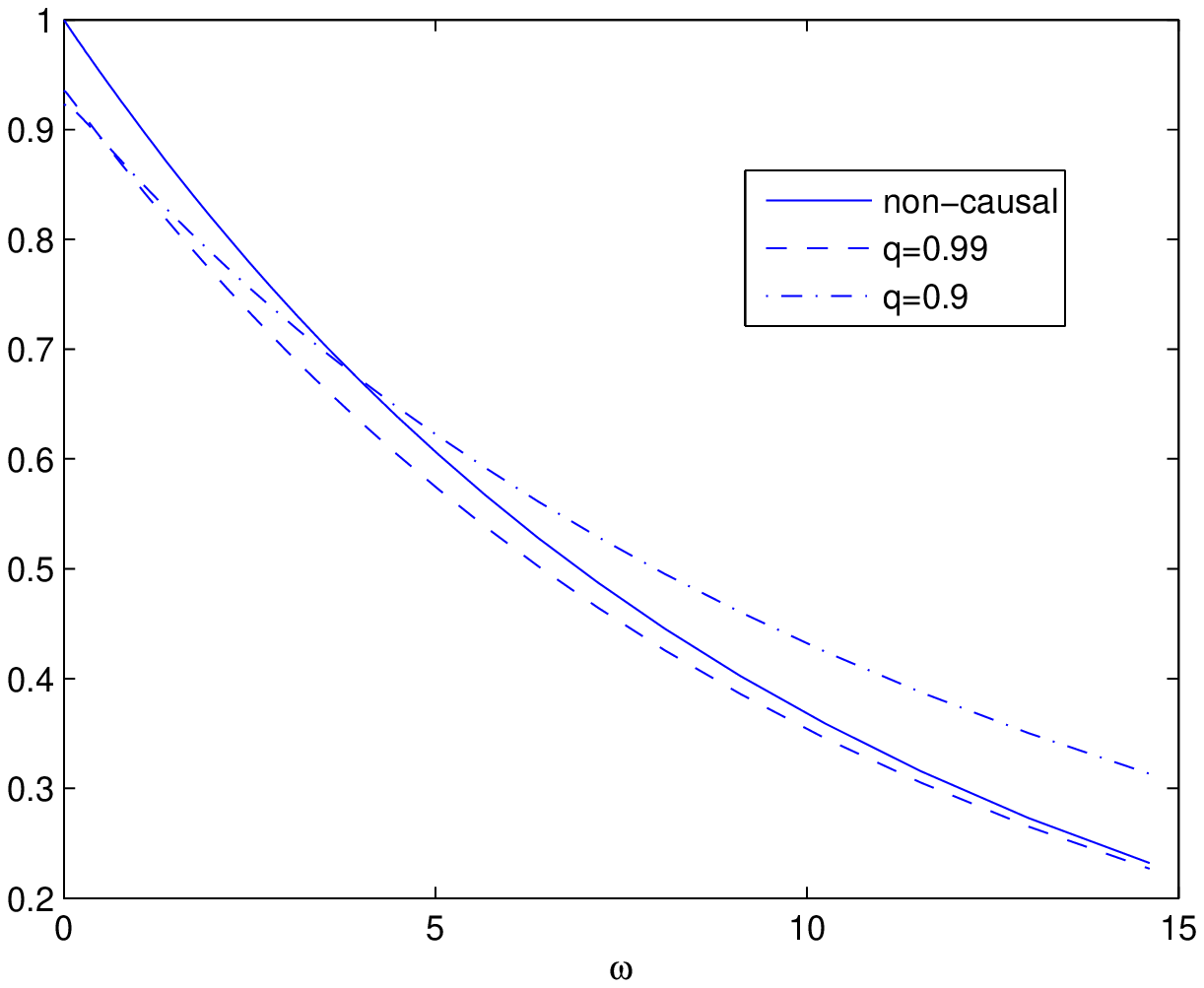,height=8.5cm}} \caption[]{\sm
Gain decay: shapes of $|M_{\mu}(i\o)|=\exp(-\mu|\o|)$ for non-causal
filter with $\mu=0.1$ and $|\HH_{\a,\b,q}(i\o)|$ for causal filters
(\ref{K}) with $q=0.99$, $b=0.01$ and $q=0.9$, $b=0.1$ respectively,
with $\a= \mu/\cos(q\pi/2)$. } \vspace{0cm}\label{fig0}\end{figure}

The sequence
$\{\HH_\lambda(\cdot)\}$ introduced in the proof above is such that
$\||H_\lambda(i\o)|-\exp(-\mu|\o|)\|_{L_2(\R)}\to 0$ as $\lambda\to
+\infty$, i.e., it approximates the gain of the non-causal smoothing
filter with the frequency response $M_{\mu}(i\o)=\exp(-\mu|\o|)$.
This sequence corresponds to  a sequence $\{\HH_{\a,\b,q}(\cdot)\}$
such that $q\to 1$, $\b\to 0$, $\a= \mu/\cos(q\pi/2)$. Figure
\ref{fig0} shows the shapes of gain curves
$|M_{\mu}(i\o)|=\exp(-\mu|\o|)$  for the reference non-causal filter
with $\mu=0.1$ and  $|H_{\a,\b,q}(i\o)|$ for sub-ideal causal
filters (\ref{K}) with $q=0.99$ and $\b=1-q=0.01$ and $q=0.9$,
$\b=1-q=0.1$ respectively. In both cases, $\a= \mu/\cos(q\pi/2)$ was
used. As expected, damping on higher frequencies is more effective
for the non-causal filter than for causal ones, and is more
effective for $q=0.99$ than for $q=0.9$. It can be illustrated as
the following: for $\o=1000$, the ratio $
\frac{|\HH_{\a,\b,q}(i\o)|}{|M_{\mu}(i\o)|}$ is found to be 1.47 and
38.65 for $q=0.99$  and $q=0.9$ respectively.

\par
\begin{figure}[ht]
\centerline{\psfig{figure=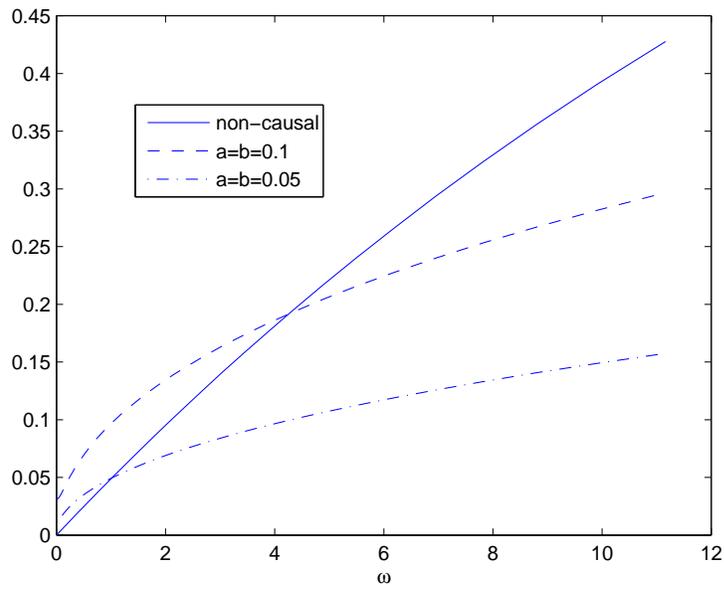,height=8.5cm}}
%\vspace{0.5cm}\centerline {\psfig{figur\href{}{}e=aminusd1.pdf,height=6.0cm}}
\vspace{0.5cm} \caption[]{\sm Approximation of identity operator:
shapes of error curves $|M_{\mu}(i\o)-1|$ and $|\HH_{a,b,q}(i\o)-1|$
respectively
 for non-causal filter with $\mu=0.05$  and for causal filters (\ref{K}) with
 $\a=\b=0.1$ and $\a=\b=0.05$, with $q=0.5$.} \vspace{0cm}
\label{fig1}
\end{figure}

 Figure \ref{fig1} illustrates Corollary \ref{corr1} and shows the shapes of error curves for approximation of
 identity operator on low frequencies. More precisely, it shows
 $|M_{\mu}(i\o)-1|=\left|e^{-\mu|\o|}-1\right|$ for the reference non-causal filter
with $\mu=0.05$ and $|\HH_{a,b,q}(i\o)-1|$  for sub-ideal causal
filters (\ref{K}) with
 $\a=\b=0.1$ and $\a=\b=0.05$ respectively, with $q=0.5$.

\begin{figure}[ht]
\centerline{\psfig{figure=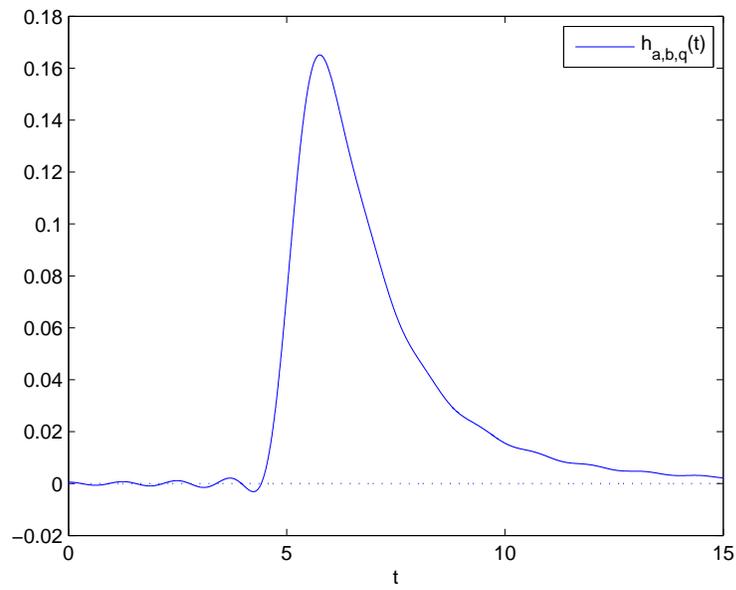,height=8.5cm}}
%\vspace{0.5cm}\centerline {\psfig{figur\href{}{}e=aminusd1.pdf,height=6.0cm}}
\vspace{0.5cm} \caption[]{\sm Impulse response
$h_{\a,b,q}(t)=(\F^{-1}H_{\a,b,q})(t)$ for causal filter (\ref{K})
with $q=0.9$, $b=0.1$, $a=1/\cos(q\pi/2)=6.3925$.} \vspace{0cm}
\label{fig3}
\end{figure}

Figure \ref{fig3} shows an example of impulse response
$h_{\a,b,q}(t)=(\F^{-1}H_{\a,b,q})(t)$ calculated  as the inverse
Fourier transform for causal filter (\ref{K}) with $q=0.9$, $b=0.1$,
$a=1/\cos(q\pi/2)=6.3925$. It can be seen that the impulse response
function almost vanishes on some interval near zero, i.e., it is
close to a causal impulse response with delay. (However, it does not
become a response with delay). There is a reason for this: if $q\to
1$ and $b\to 0$ then, for a given $\a,c>0$, $H_{\a,b,q}(p)\to e^{\a
p}$ uniformly in the domain $\{p\in \C^+: |p|\le c\}$.

\par
It can be noted that the phase shift for the frequency response
function is large for large $|\o|$, and it is increasing  when $q\to
1$. This does not affect much the performance of the filter since
the gain is small for these large $|\o|$ and Condition (a1) is
ensured.

%[a]=filter2011qq(0.15,0.7,0.9,160)
\section{Conclusion}
The paper proposes a family of causal smoothing filters with almost
exponential damping of the energy on the higher frequencies and with
the frequency response that can be selected to be arbitrarily close
to the real unity uniformly on an arbitrarily large interval. These
filters are sub-ideal meaning that a faster decay of the frequency
response would lead to the loss of causality; this is because they
approximate non-causal filters with exponential rate of decay.  A
possible application is in interpolation and forecast algorithms.
The transfer functions obtained are not rational functions; it would
be interesting to consider their approximation by the rational
functions. Another problem is the transition to discrete time
processes. We leave it for future research.

%\newpage

\end{document}